# Detection and discrimination of the periodicity of prime numbers by discrete Fourier transform

## – Symphony of primes –

Levente Csoka

University of West Hungary, Institute of Wood Based Products and Technologies, 9400 Sopron, Hungary

Address correspondence to this author at the Inst. of Wood Based Products and Technologies, University of West Hungary, Hungary; Tel.: +36-99518-305, Fax: +36-99-518-302; E-mail: levente.csoka@skk.nyme.hu

ABSTRACT. A novel representation of a quasi-periodic modified von Mangoldt function $L(n)$ on prime numbers and its decomposition into Fourier series has been investigated. We focus on some particular quantities characterizing the modified von Mangoldt function. The results indicate that prime number progression can be decomposed into periodic sequences. The main approach is to decompose it into sin or cosine function. Basically, it is applied to extract hidden periodicities in seemingly quasi periodic prime function. Numerical evidences were provided to confirm the periodic distribution of primes.

I. The von Mangoldt function

The von Mangoldt arithmetic function is defined by

$$\Lambda(n) \equiv \begin{cases} \log p & \text{if } n = p^k \text{ for } p \text{ a prime and integer } k \geq 1 \\ 0 & \text{otherwise.} \end{cases} \qquad 1$$

This function carries the information about the properties of the primes and approximate their weighted characteristic function. The discrete Fourier transform of the von Mangoldt function is defined by the Riemann-Weil explicit formula

$$\sum_\gamma h(\gamma) = h\left(\frac{i}{2}\right) + h\left(-\frac{i}{2}\right) -$$

$$-g(0)\log \pi + \frac{1}{2\pi}\int_{-\infty}^{\infty} h(r)\frac{\Gamma'}{\Gamma}\left(\frac{1}{4} + \frac{1}{2}ir\right)dr - 2\sum_{n=1}^{\infty} \frac{\Lambda(n)}{\sqrt{n}} g(\log n) \qquad 2$$

and gives a spectrum with spikes at ordinates equal to the imaginary part of the Riemann zeta function zeros (Hejhal 1976).

II. The modified von Mangoldt function

Starting with the function above, let us consider a special, simple case of von Mangoldt function at for $k = 1$,

$$L(n) \equiv \begin{cases} \ln(e) & if\ n = p\ for\ p\ a\ prime \\ 0 & otherwise. \end{cases} \qquad 3$$

The function $L(n)$ is a representative way to introduce the set of primes with weight of unity attached to the location of a prime. From the summation of that function we can see that primes contribute equally and that representation exist. The sum of the modified von Mangoldt function $L(n)$ is same as the Riemann hypothesis (Devenport 1980, Vardi 1991) with a slight error:

$$\psi(n) = \int_2^{p_n} \frac{1}{\ln t} dt = n + O(\sqrt{n}\,(\ln x)^2) \equiv \sum_{n \to \infty} L(n) \equiv p_n \qquad 4$$

but the $L(n)$ function predict the number of primes with no error, e.g. the primes are not randomly distributed (large numbers' law). Therefore, this function is oscillatory but not diverging. The modified von Mangoldt function $L(n)$ is related to Dirichlet seriesor the Riemann zeta function by

$$\zeta(a) = \sum_{n=1}^{\infty} \frac{1}{n^a} = \sum_{n=1}^{\infty} \frac{L(n)}{n^a} \quad \Re(a) > 1. \qquad 5$$

The discrete Fourier transform of the modified von Mangoldt function $L(n)$ gives a spectrum with periodic, real parts as spikes at the frequency axis ordinates. Hence the modified von Mangoldt function $L(n)$ and the prime series can be approximated by periodic sin or cos functions. Let us consider that natural numbers form a discrete function $g(n) \to g(n_k)$ and $g_k \equiv g(n_k)$ by uniform sampling at every sampling point $n_k \equiv k\Delta$, where $\Delta$ is the sampling period and $k = 2, 3, .., L - 1$. If the subset of discrete prime numbers is generated by the uniform sampling of a wider discrete set, it is possible to obtain the Discrete Fourier Transform (DFT) of $L(n)$. The sequence of $L(n)$ prime function is transformed into another sequence of N complex numbers according to the following equation:

$$F\{L(n)\} = X(\nu) \qquad 6$$

where $L(n)$ is the modified von Mangoldt function. The the operator $F$ is defined as:

$$X(\nu) = \sum_{n=2}^{N-1} L(k\Delta) e^{-i(2\pi\nu)(k\Delta)}, \qquad 7$$

where $\nu = l\,\Delta f$, $l = 0, 1, 2, \ldots, N - 1$, and $\Delta f = 1/L$, $L$=length of the modified von Mangoldt function.

The amplitude spectrum of equation 7, describes completely the discrete-time Fourier transform of the $L(n)$ prime function of quasi periodic sequences, which comprises only discrete frequency component. The important properties of the amplitude spectrum are as follows:

1) the ratio of consecutive frequencies $(f_t, f_{t+1})$ and maximum frequency $(f_{max})$ converges $\lim_{t \to f_s/2} \frac{f_{t+1}}{f_t} = \frac{f_t + (1/L)}{f_t} = \frac{f_{max}}{f_t} = 1$, where $f_s = \frac{1}{\Delta}$ represents the positive half-length of amplitude spectrum. This is similar to the prime number theorem: $\lim_{x \to \infty} \frac{\pi(x) \ln x}{x} = 1.$      8

2) the reciprocate of the frequencies $(f_t)$ converges $\lim_{t \to f_s/2} \frac{1}{f_t} = \frac{1}{f_s}$     9

3) the consecutive frequencies on the DFT amplitude spectrum describe a parabolic or Fermat's spiral (Fig.1.), starting from 0 with increasing modulus (polar form $r = af_t$), where $a = 1$, independently from the length of the prime function and its parametric equations can be written as: $x(f_t) = f_t \cos(f_t), y(f_t) = f_t \sin(f_t)$     10

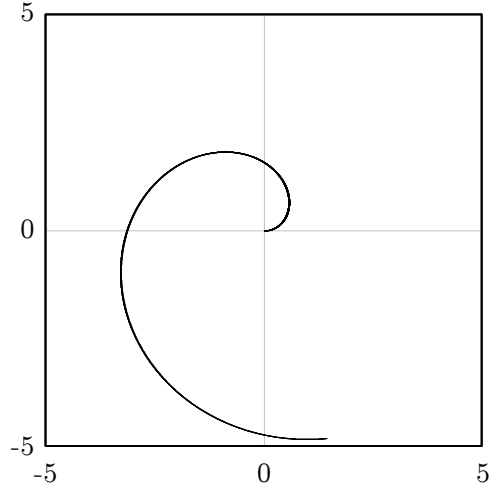

Fig. 1. Parabolic or Fermat spiral using the consecutive frequencies. The arc length of the spiral is a function of the sampling interval.

4) the DFT spectrum of the primes exhibits highly ordered and symmetrical frequency distribution (Fig.3.)
5) the DFT of the modified von Mangoldt function will be periodic (Fig.2.) independently the length of integer sequences used ($v \geq N, v = 0,1, \ldots, N-1$), $X(v) = X(v + N) = X(v + 2N), \ldots, = X(v + zN)$

$$X(zN + v) = \sum_{n=0}^{N-1} L(n) e^{-i(2\pi nv)(zN+v)} = \sum_{n=0}^{N-1} L(n) e^{-i(2\pi nv)} e^{-i2\pi zn} \quad 11$$

$z$ and $n$ are integers and $e^{-i2\pi zn} = 1$.

$$X(v + zN) = \sum_{n=0}^{N-1} L(n) e^{-i(2\pi nv)} = X(v) \quad 12$$

III.    Spectral evidence

For completeness we give the compelling spectral evidence for the truth that Fourier decomposition of the modified von Mangoldt function is periodic. Fourier transformation has been used for studying the behaviour of prime numbers for a long time and it was employed here to reveal the periodic distribution of the primes. Figure 2 shows that the corresponding amplitude spectrum of the modified von Mangoldt function exhibits a well-defined periodic signature. This suggests that the series of prime numbers is not an arbitrary function. Creating the superposition from the Fermat spiral frequencies by adding their amplitude and phase together, the original prime function can be restored.

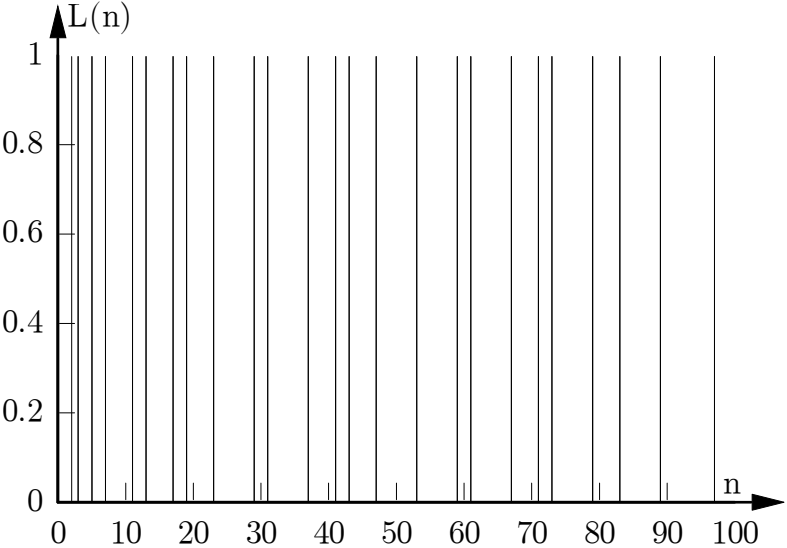

Fig. 2. Representation of the modified von Mangoldt function up to 100. Spikes appear at the location of prime numbers.

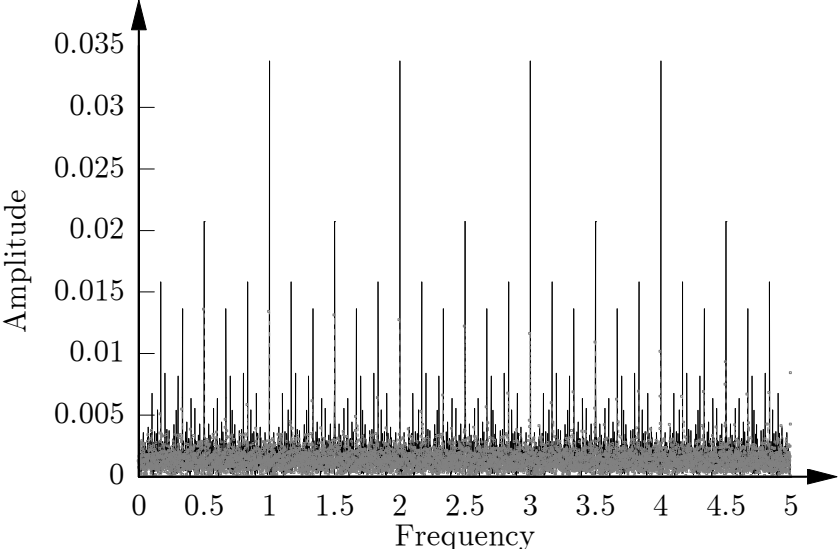

Fig. 3. The Discrete Fourier amplitude spectrum of the modified von Mangoldt function a) DFT of primes between 2-10.000 (continuous black line) and b) DFT of primes between 1.0200*10^7-1.0201*10^7 (dotted grey line)

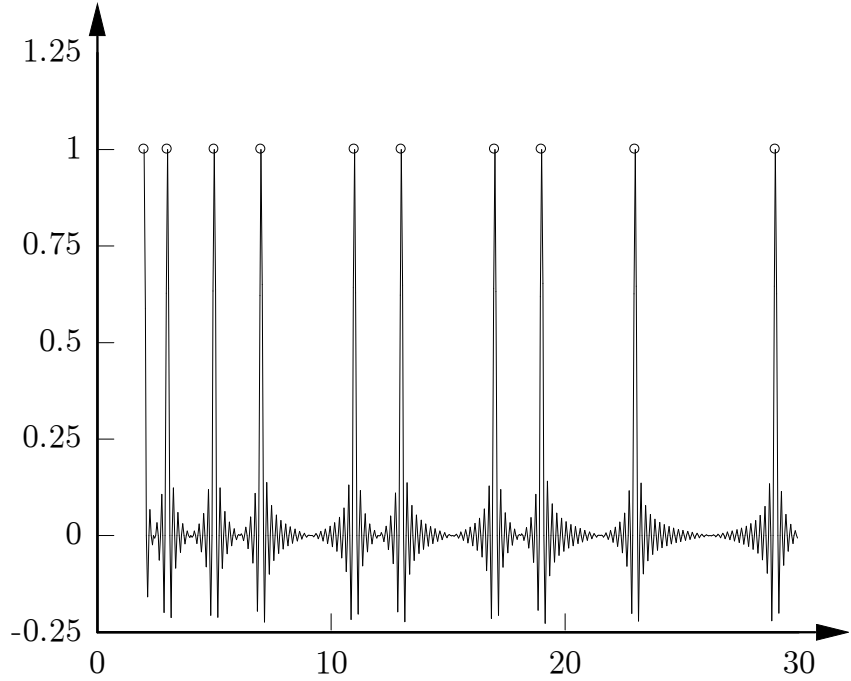

Fig. 4. Reconstruction of the modified von Mangoldt function from the Fourier decomposed periodic waves between 2-10.000 (magnified from 2-30). The circles at the top of the spikes indicates the location of prime numbers.

Creating the superposition from the Fermat spiral frequencies by adding their amplitude and phase together, the original prime function can be arising.

$$L(n) = \sum_{f_t=0}^{f_{max}} A_t \sin(f_t + \omega n) \qquad 13$$

From the wave sequence above, it can be seen that all primes have a last digit of 1, 3, 5 and 7 or 9 [4], so the music of the primes billow over the set of natural numbers and able to show some arithmetical progression, pattern of prime numbers [5]. Figure 5 shows that new prime numbers are decreasing to infinity and that how many new primes we have for example in a 1000 number around 1000, 10^3, 10^4 etc., so less and less amount of new frequencies requires to describe their pattern.

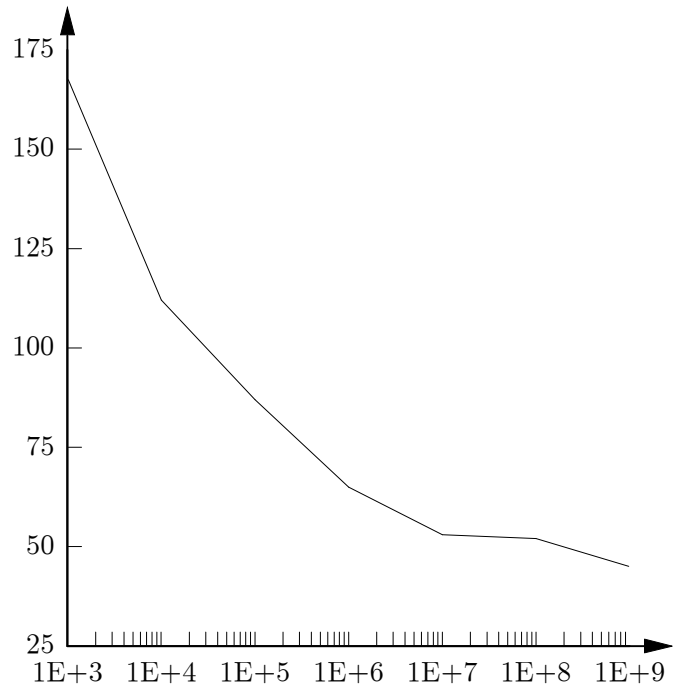

Fig. 5. Number of new primes in a 1000 number interval.

Conclusion

We have presented a novel method for the detection of prime repetitions on the set of natural numbers. The method is based on a DFT of the modified von Mangoldt function. The relevance of the method for determination of the conformational sequences is described.